\documentclass[a4paper,12pt]{article}
\usepackage[latin1]{inputenc}
\usepackage{amssymb}
\title{{\bf  Local sentences and Mahlo cardinals } }
\author{Olivier Finkel$^1$    ~~~~ and  ~~~~  Stevo Todorcevic$^{2, 3}$
\\ 
\\
$^1${\it Equipe Modèles de Calcul et Complexité}  
 \\ {\it Laboratoire de l'Informatique du Parallélisme}\footnote{UMR 5668 - CNRS - ENS Lyon - UCB Lyon - INRIA }
 \\  CNRS et Ecole Normale Supérieure de Lyon
 \\ 46, Allée d'Italie 69364 Lyon Cedex 07, France.
\\ Olivier.Finkel@ens-lyon.fr 
\\
\\  $^2${\it Equipe de Logique Math\'ematique }
 \\ U.F.R. de Math\'ematiques, Universit\'e Paris 7 \\  2 Place Jussieu 75251 Paris
 cedex 05, France.
\\ stevo@logique.jussieu.fr
\\
\\ $^3${\it  Department of Mathematics} 
\\ University of Toronto, Toronto, Canada
M5S 2E4.
}

\date{}
\begin{document}

\newtheorem{The}{Theorem}[section]
\newtheorem{Pro}[The]{Proposition}
\newtheorem{Deff}[The]{Definition}
\newtheorem{Lem}[The]{Lemma}
\newtheorem{Rem}[The]{Remark}
\newtheorem{Exa}[The]{Example}
\newtheorem{Cor}[The]{Corollary}

\newcommand{\vp}{\varphi}
\newcommand{\lb}{\linebreak}
\newcommand{\fa}{\forall}
\newcommand{\Ga}{\Gamma}
\newcommand{\Gas}{\Gamma^\star}
\newcommand{\Gao}{\Gamma^\omega}
\newcommand{\Si}{\Sigma}
\newcommand{\Sis}{\Sigma^\star}
\newcommand{\Sio}{\Sigma^\omega}
\newcommand{\ra}{\rightarrow}
\newcommand{\hs}{\hspace{12mm}

\noi}
\newcommand{\lra}{\leftrightarrow}
\newcommand{\la}{language}
\newcommand{\ite}{\item}
\newcommand{\Lp}{L(\varphi)}
\newcommand{\abs}{\{a, b\}^\star}
\newcommand{\abcs}{\{a, b, c \}^\star}
\newcommand{\ol}{ $\omega$-language}
\newcommand{\orl}{ $\omega$-regular language}
\newcommand{\om}{\omega}
\newcommand{\nl}{\newline}
\newcommand{\noi}{\noindent}
\newcommand{\tla}{\twoheadleftarrow}
\newcommand{\de}{deterministic }
\newcommand{\proo}{\noi {\bf Proof.} }
\newcommand {\ep}{\hfill $\square$}

\maketitle

\begin{abstract}
\noi    Local sentences were introduced by Ressayre in \cite{ress}
who proved certain remarkable stretching
theorems establishing the equivalence between the existence of
finite models for these sentences and the existence of some infinite
well ordered models.  Two of these stretching theorems were only
proved under certain large cardinal axioms but the
question of their exact (consistency) strength was left open in
\cite{fr}. Here, we solve this problem, using a combinatorial result
of J. H. Schmerl \cite{sch}. In fact, we show that the
stretching principles are equivalent to the existence of
$n$-Mahlo cardinals for appropriate integers $n$. This is done
by proving first that for each integer $n$, there is a local
sentence $\phi_n$ which has  well ordered models of order type
$\tau$, for every infinite ordinal $\tau > \om$ which is not an
$n$-Mahlo cardinal.
\end{abstract}

\noi {\small {\bf Keywords:}     Local sentences; stretching
principles; well ordered models;  Mahlo cardinals; spectra of
sentences.     }

\section{Introduction}

Local  sentences were introduced
by  J.-P. Ressayre who proved
 some remarkable links between the finite and the infinite model theory of
these sentences,   \cite{ress}.  A local sentence is a first order sentence
which is equivalent to a universal sentence and satisfies some semantic
restrictions:
closure in its models takes a finite number of steps.
Assuming that a binary relation symbol belongs to the signature of a local
sentence $\vp$
and is interpreted by a linear order in every model of $\vp$, the
stretching theorems state that the existence of certain well ordered models of
$\vp$ is equivalent to
the existence of a finite model of $\vp$, generated by
some particular kind of indiscernibles, like special, remarkable, semi-monotonic or
monotonic ones (see \cite{fr} for a precise definition).
Two of these stretching theorems were
only  proved under large cardinal axioms  but the question of  their exact (consistency) strength was left open in \cite{fr}.
We solve  here this problem, using a combinatorial result of J. H. Schmerl which characterizes 
 $n$-Mahlo cardinals in \cite{sch}. We show that the two stretching theorems are
in fact equivalent to the existence, for each integer $n$, of an $n$-Mahlo cardinal.
This is done by proving first that for each integer $n$, 
there is a local sentence $\phi_n$ which has some well ordered models
of order type $\tau$, for every infinite ordinal $\tau > \om$ which is not an $n$-Mahlo cardinal.
Using this result we show also that a kind of hanf number 
$\mu$ for local sentences defined by Ressayre in \cite{ress} 
is in fact equal to the first $\om$-Mahlo cardinal,
if such a cardinal exists.

\hs The paper is organized as follows. In section 2 we recall some definitions and stretching theorems. In section 3 we prove the existence of a local
sentence $\phi_n$ which has some well ordered models
of order type $\tau$, for every infinite ordinal $\tau > \om$ which is not an $n$-Mahlo cardinal. In section 4 we solve the question of the
 exact (consistency) strength of some stretching principles for local sentences.

\section{Stretching theorems for local sentences}

In this paper the (first order) signatures are finite, always contain one
binary
predicate symbol $=$ for equality, and can contain both functional
and relational symbols.

\noi When $M$ is a structure in a signature $\Lambda$, $|M|$ is the domain of $M$.
\nl If $S$ is  a function, relation  or  constant symbol in $\Lambda$, then
 $S^M$ is the interpretation in the structure
$M$ of  $S$.
\nl Notice that, when the meaning is clear, the superscript $M$ in
$S^M$ will be sometimes omitted in order to simplify the presentation.

 For  $X\subseteq |M|$,
we define:
\nl $cl^1(X, M)=X \cup \bigcup_{\{f ~{\rm k-ary~ function~ of ~} \Lambda
~\} } ~f^M(X^k)
\cup \bigcup_{\{a ~{\rm ~ constant~ of ~} \Lambda ~\} } a^M $
\nl $cl^{n+1}(X, M)=cl^1(cl^n(X, M), M) \quad {\rm ~ for ~an ~integer~}
n\geq 1$
\nl and $cl(X, M)=\bigcup_{n\geq 1} cl^n(X, M)$ is the closure of $X$ in
$M$.

 The
signature of a first order sentence $\varphi$, i.e. the set of
non logical symbols appearing in $\varphi$, is denoted S($\varphi$).

\begin{Deff}\label{defloc} A first order sentence $\varphi$
is local if and only if:
\begin{enumerate}
\item[(a)] $M\models\varphi$ and $X\subseteq|M|$ imply
$cl(X,M)\models\varphi$
\item[(b)] $\exists n\in \mathbb{N}$ such that $\forall M$, if $ M\models
\varphi$ and
$X \subseteq |M|$,
then $ cl(X,M)=cl^n(X, M)$, (closure in models of $\varphi$ takes at most
$n$ steps).

\end{enumerate}
\end{Deff}

\noi For a local sentence $\varphi$, $n_\varphi$ is the
smallest integer $n\geq 1$ satisfying $(b)$ of the above definition.
In this definition, $(a)$ implies that a local sentence
$\varphi$ is always equivalent to a universal
sentence, so we may assume that this is always the case.

\begin{Exa}
Let $\varphi$ be the sentence, (already given in \cite{wollic}),  in the signature $S(\varphi)=\{<,
P, i, a\}$,
where $<$ is a binary relation symbol, $P$ is a unary relation symbol, $i$
is a unary
function symbol, and $a$ is a constant symbol, which is the conjunction of:

\begin{enumerate}
\ite[(1)] $\fa xyz [ (x\leq y \vee y\leq x) \wedge ((x\leq y \wedge y\leq
x) \lra x=y) \wedge
((x\leq y \wedge y\leq z) \ra x\leq z) ]$,
\ite[(2)] $\fa x y [ ( P(x) \wedge \neg P(y) ) \ra x<y ]$,
\ite[(3)] $\fa x y [ ( P(x) \ra i(x)=x ) \wedge ( \neg P(y) \ra
P(i(y)) ) ]$,
\ite[(4)] $\fa x y [ ( \neg P(x) \wedge \neg P(y) \wedge x \neq y ) \ra
i(x) \neq i(y) ]$,
\ite[(5)] $\neg P(a)$.
\end{enumerate}
\end{Exa}

\noi We now explain the meaning of the above sentences $(1)$-$(5)$.
\nl Assume that $M$ is a model of $\vp$.
The sentence $(1)$ expresses that $<$ is interpreted in $M$ by a linear
order; $(2)$ expresses
that $P^M$ is an initial segment of the model $M$; $(3)$ expresses that the
function $i^M$
is trivially defined by $i^M(x)=x$ on $P^M$ and is defined from $\neg P^M$
into $P^M$.
$(4)$ says that $i^M$ is an injection from $\neg P^M$ into $P^M$ and $(5)$
ensures that
the element $a^M$ is in $\neg P^M$.
\nl The sentence $\varphi$ is a conjunction of universal sentences thus it
is
equivalent to a universal one, and closure in its models takes at most
two steps (one adds the constant $a$ in the first step then takes the closure
under the function $i$).
Thus $\varphi$ is a local sentence.
\nl If we consider only the order types of {\it well ordered} models of
$\varphi$,
we can easily see that $\varphi$ has a model of order type $\alpha$, for
every finite ordinal
$\alpha \geq 2$ and for every infinite ordinal $\alpha$ which is not a cardinal.

\hs The reader may
also find many other examples of local
sentences
 in the papers \cite{ress,fr,loc,cloloc,wollic}.  Notice that local sentences play a role in defining many classes of 
formal languages of finite or infinite words which are important in theoretical computer science, like the classes of regular or quasirational 
languages; the latter one  forms of rich subclass of the class of  context-free languages, see \cite{ress,loc,cloloc}. 

\hs  From now on we shall assume that the signature of local sentences
contain
a binary predicate $<$ which is interpreted by a linear ordering in all
of their models.
\nl We recall  the stretching theorem for local sentences.
Below, semi-monotonic, special, remarkable,
and monotonic indiscernibles are particular kinds of indiscernibles
which are precisely defined in \cite{fr}.

\begin{The}[\cite{fr}]\label{stretching}
For each local sentence $\varphi$ there exists a positive
integer $N_\varphi$ such that
\begin{enumerate}
\ite[(A)] $\varphi$ has arbitrarily large finite models if and only if
$\varphi$ has
an infinite model if and only if $\varphi$ has a finite model generated by
$N_\varphi$
indiscernibles.
\ite[(B)] $\varphi$ has an infinite well ordered model
if and only if $\varphi$ has a finite model generated by $N_\varphi$
semi-monotonic
indiscernibles.
\ite[(C)] $\varphi$ has a model of order type $\om$
if and only if $\varphi$ has a finite model generated by $N_\varphi$ special
indiscernibles.
\ite[(D)] $\varphi$ has well ordered models of unbounded order types in
the ordinals
if and only if $\varphi$ has a finite model generated by $N_\varphi$
monotonic
indiscernibles.
\ite[(E)] $\varphi$ has for every infinite cardinal $k$ a model of order type $k$ if and only if
$\varphi$ has a finite model generated by $N_\varphi$
monotonic and remarkable
indiscernibles.
\ite[(E')] $\varphi$ has for every infinite cardinal $k$ a model of order type $k$  in which a distinguished predicate
$P$ is cofinal if and only if
$\varphi$ has a finite model generated by $N_\varphi$
monotonic and special
indiscernibles belonging to $P$.

\end{enumerate}
\noi To every local sentence $\varphi$ and every ordinal $\alpha$ such that
$\om \leq \alpha <\om^\om$ one can associate by an effective procedure a
local sentence
$\vp_\alpha$, a unary predicate symbol $P$ being in the signature
S($\varphi_\alpha$),
such that:
\begin{enumerate}
\ite[($C_\alpha$)] $\varphi$ has a well ordered model of order type
$\alpha$
if and only if $\varphi_\alpha$ has a finite model $M$ generated by
$N_{\varphi_\alpha}$
semi-monotonic indiscernibles  belonging to $P^M$.
\end{enumerate}

\end{The}

\noi It is proved in \cite{fr} that the integer $N_\varphi$ can be
effectively computed from $n_\varphi$ and $q$ where $$\varphi=\fa
x_1 \ldots \fa x_q \theta(x_1, \ldots, x_q)$$ and $\theta$ is a quantifier free  formula. 
If $v(\varphi)$ is the maximum number of variables of terms of complexity $\leq
n_\varphi +1$ (resulting by at most $n_\varphi +1$ applications of
function symbols) and $v'(\varphi)$ is the maximum number of
variables of an atomic formula involving terms of complexity $\leq
n_\varphi +1$ then 
$$N_\varphi = max \{ 3v(\varphi) ; v'(\varphi) + v(\varphi) ;
q.v'(\varphi) \}.$$

\hs It is also proved in \cite{fr} that actually some equivalences
of this stretching theorem hold only under strong axioms of
infinity:

\hs $(E)$ is provable in {\bf ZF} + existence for each integer $n$ of an $n$-Mahlo cardinal; but not in
{\bf ZF}, for it implies the existence of an inaccessible cardinal.

\hs $(E')$ is provable in {\bf ZF} + the scheme asserting for each standard integer
 $n$ the existence of an $n$-Mahlo cardinal; and $(E')$ implies the consistency of this scheme.

\hs The question of the exact (consistency) strength of $(E)$ and  $(E')$ was left open in \cite{fr} and is
solved in this paper.

\section{Infinite spectra of local sentences}

\noi We are mainly interested in this paper by {\it well ordered} models of local sentences,
so we now recall the notion of spectrum of a local sentence $\vp$.
As usual the class of all ordinals is denoted by {\bf On}.

\begin{Deff}
Let $\vp$ be a local sentence; the spectrum of $\vp$ is
$$Sp(\vp) = \{ \tau \in {\bf On} \mid \vp \mbox{ has a model of order type } \tau \}$$
\noi and the infinite spectrum of $\vp$ is
$$Sp_\infty(\vp) = \{ \tau \in {\bf On} \mid \tau \geq \om \mbox{ and }
\vp \mbox{ has a model of order type } \tau \}$$
\end{Deff}

The following result was proved in \cite{fr}:

\begin{The}\label{phi0}
There exists a local sentence $\phi_0$ such that
$$Sp_\infty(\phi_0) = \{ \tau > \om \mid \tau \mbox{ is not an inaccessible cardinal } \}.$$
\end{The}

\noi We are going firstly to extend this result by proving the following one.

\begin{The}\label{phin}
For each integer $n\geq 1$, there exists a local sentence $\phi_n$ such that
$$Sp_\infty(\phi_n) = \{ \tau > \om \mid \tau \mbox{ is not an  } \mbox{n-Mahlo cardinal } \}.$$
\end{The}

\noi In order to construct the local sentences $\phi_n$, we shall use a combinatorial result of
J. H. Schmerl which gives a characterization of $n$-Mahlo cardinals \cite{sch}. Notice that in
\cite{sch}, the now usually called $n$-Mahlo cardinals are just called n-inaccessible.

\hs As usual, the set of subsets of cardinality $n$ of a set $X$ is denoted by $[X]^n$.
If $C$ is a partition of $[X]^n$ then $Y \subseteq X$ is $C$-homogeneous iff every two elements
of  $[Y]^n$ are $C$-equivalent, i.e. are in the same set of $C$.

\begin{Deff} For an integer $n\geq 1$ and an ordinal $\alpha$, let  $P(n, \alpha)$ be the class of
infinite cardinals $k$ which have the following property: Suppose for each $\nu < k$ that $C_\nu$ is a partition
of $[k]^n$ and card($C_\nu$) $ <k$ then there is $X \subseteq k$ of length $\alpha$ such that for each
 $\nu \in X$, the set $X - (\nu + 1)$ is $C_\nu$-homogeneous.
\end{Deff}

\begin{The}[\cite{sch}]
Let $n\geq 1$ be an integer and $k$ be an infinite cardinal. Then $k \in P(n+2, n+5)$ if and only if
$k=\om$ or $k$ is an $n$-Mahlo cardinal.
\end{The}

\noi  We are going to construct a local sentence $\theta_n$ such
that, for each regular infinite cardinal $\kappa$, it holds that:

( $\theta_n$ has a model of order type $\kappa$ ) if and only if (
$\kappa \notin P(n+2, n+5) )$

\hs So we have to express that, for each $\nu < \kappa$,  there is a
partition $C_\nu$ of $[\kappa]^{n+2}$ with
 card($C_\nu$) $ <\kappa$, such that,  for all subsets $X$ of $\kappa$ having $n+5$ elements,  there exists
 $\nu \in X$ such that  $X - (\nu + 1)$ is not $C_\nu$-homogeneous.

\hs We shall firstly express by the following sentence $\psi_1$ that a model $M$ is divided into successive segments.
The signature of  $\psi_1$ is $\{<, I, P\}$, where $I$ is a unary function and $P$ is a unary predicate.
 $\psi_1$ is the  conjunction of :

\begin{itemize}
\ite[(1)] ( $<$ is a linear order ),

\ite[(2)]  $\fa yz [  y \leq I(y) $ and $ (y \leq z \ra I(y) \leq I(z)) $ and
$ ( y  \leq z \leq I(y)  \ra I(z)=I(y) ) ] $,

\ite[(3)]   $\fa y [  I(y)=y  \leftrightarrow P(y) ]$.

\end{itemize}

\noi  In a model $M$ of  $\psi_1$,
the function $I^M$ is constant on each of these segments and the image $I^M(x)$
of an element
$x$ is the last  element of the segment containing $x$.
We have added that  the set of the last elements of the successive segments is the subset $P^M$ of $|M|$.
\nl The sentence $\psi_1$ is equivalent to a universal sentence and closure (under the function $I$)
 in its models takes at most one step thus $\psi_1$ is a local sentence.

\hs If $M$ is a well ordered model of $\psi_1$ whose order type is a
regular cardinal $\kappa$, then the set $P^M$ is cofinal in $\kappa$
so the order type of $(P^M, <^M)$ will be also $\kappa$. The set
$P^M$ will then be identified with $\kappa$ and each segment defined
by the sentence $\psi_1$ will be of cardinal smaller than $\kappa$
(because it is bounded in the model $M$ by the last element of the
segment).

\hs We are now going to express, using a $(n+3)$-ary function $f$,
that, for each $\nu < \kappa$, there is a partition $C_\nu$ of
$[\kappa]^{n+2}$ with  card($C_\nu$) $ <\kappa$. We shall use the
following sentence $\psi_2$ in the signature $S(\psi_2)=\{<, I, P,
f\}$, which is the conjunction of :

\begin{itemize}

\ite[(1)]  $\fa \nu y_1 y_2 \ldots y_{n+2}
[ f(\nu, y_1,  y_2, \ldots, y_{n+2})=f(I(\nu), I(y_1),  I(y_2), \ldots, I(y_{n+2}))]$,

\ite[(2)]  $\fa \nu y_1 y_2 \ldots y_{n+2} [ (  \bigvee_{1 \leq i < j \leq n+2} \neg (I(y_i) < I(y_j))  )   \ra
 f(\nu, y_1,  y_2, \ldots, y_{n+2})=I(\nu) ]$,

\ite[(3)]  $\fa \nu y_1 y_2 \ldots y_{n+2} [ (  \bigwedge_{1 \leq i < j \leq n+2}  (I(y_i) < I(y_j))   )   \ra
I( f(\nu, y_1,  y_2, \ldots, y_{n+2}) ) =I(\nu) ]$,

\ite[(4)]  $\fa \nu y_1 y_2 \ldots y_{n+2} [ (  \bigwedge_{1 \leq i < j \leq n+2}  (I(y_i) < I(y_j))   )  \ra
\neg P( f(\nu, y_1,  y_2, \ldots, y_{n+2}) )  ]$.
\end{itemize}

\noi  If $M$ is a well ordered model of $\psi_1\wedge \psi_2$ whose
order type is a regular cardinal $\kappa$, then we have, for each
$\nu < \kappa$ (identified with  $\nu \in P^M$ so $\nu=I(\nu)$)  a
partition $C_\nu$ of $[\kappa]^{n+2}$ which is defined by the
function $f$. For all elements $y_1 < y_2 < \ldots < y_{n+2}$ and
$y'_1 < y'_2 < \ldots < y'_{n+2}$ of $\kappa$ (so all elements $y_i$
and $y'_i$ are in $P^M$) the two sets $\{y_1, y_2, \ldots y_{n+2}\}$
and $\{y'_1, y'_2, \ldots y'_{n+2}\}$     are         in the same
set of $C_\nu$ iff
$$f(\nu, y_1,  y_2, \ldots, y_{n+2}) =  f(\nu, y'_1,  y'_2, \ldots, y'_{n+2}).$$
\noi But the elements $f(\nu, y_1,  y_2, \ldots, y_{n+2})$ and
$f(\nu, y'_1,  y'_2, \ldots, y'_{n+2})$ are in the segment of $M$
whose last element is $I(\nu)$ (and they are different from
$I(\nu)$) thus it will hold that card($C_\nu$) $ <\kappa$ because we
have seen that each segment defined by $\psi_1$ will have
cardinality smaller than $\kappa$.

\hs Notice that items $(1)$ and $(2)$
above ensure that the function $f$ is always defined and that closure in models of $\psi_1\wedge \psi_2$
takes at most two steps : one takes closure under the function $I$ in one step, then closure under the
function $f$ in a second step. So the sentence  $\psi_1\wedge \psi_2$ is local.

\hs We have now to express that,  for all subsets $X$ of $\kappa$
having $n+5$ elements,  there exists
 $\nu \in X$ such that  $X - (\nu + 1)$ is not $C_\nu$-homogeneous.

\hs We shall use the following sentence $\psi_3$ in the same signature $\{<, I, P, f\}$, which is equal to :

\hs   $\fa x_1 \ldots x_{n+5} \in P [ \bigwedge_{1 \leq i < j \leq n+5} x_i < x_j  \ra \phi( x_1 \ldots x_{n+5})]$,

\hs where $\phi( x_1 \ldots x_{n+5})$ is the sentence :

\hs $ \bigvee_{\nu \in \{x_1, x_2, x_3\}}   \bigvee_{ \{y_1, \ldots, y_{n+2}\} \in [\{x_j \mid 1\leq j \leq n+5\}]^{n+2}}
 \bigvee_{ \{y'_1, \ldots, y'_{n+2}\} \in [\{x_j \mid 1\leq j \leq n+5\}]^{n+2}}
\nl [ \nu < y_1< \ldots < y_{n+2} \bigwedge \nu < y'_1< \ldots < y'_{n+2} \bigwedge
 f(\nu, y_1,  y_2, \ldots, y_{n+2})  \neq  f(\nu, y'_1,  y'_2, \ldots, y'_{n+2})  ]$.

\hs Consider now the sentence $\theta_n = \bigwedge_{1\leq i \leq 3}
\psi_i$. This sentence is local. Moreover if $\kappa$ is an infinite
regular cardinal,  then $\theta_n$ has a model of order type
$\kappa$ iff for each $\nu < \kappa$, there is  a partition $C_\nu$
of $[\kappa]^{n+2}$ with
 card($C_\nu$) $ <\kappa$, such that for all subsets $X$ of $\kappa$ having $n+5$ elements,  there exists
 $\nu \in X$ such that  $X - (\nu + 1)$ is not $C_\nu$-homogeneous.

\hs Recall now that from two local sentences $\varphi_1$ and $\varphi_2$ we can construct another local sentence
$\varphi_1 \cup \varphi_2$ such that $Sp(\varphi_1 \cup \varphi_2)=Sp(\varphi_1) \cup Sp(\varphi_2)$, \cite{fr}.

\hs Consider  now   the local sentence $\phi_0$
given in Theorem \ref{phi0} such that
$$Sp_\infty(\phi_0) = \{ \tau > \om \mid \tau \mbox{ is not an inaccessible cardinal } \}.$$
\noi The local sentence $\phi_n=\theta_n \cup \phi_0$ is a local sentence and
$$Sp_\infty(\phi_n) = \{ \tau > \om \mid \tau \mbox{ is not an  } n \mbox{-Mahlo cardinal } \}.$$
\noi This ends the proof of Theorem \ref{phin}.  \ep

\hs We can now determine the kind of Hanf number $\mu$ for local sentences defined in \cite{fr}.

\begin{Deff}
The ordinal $\mu$ is the smallest ordinal such that, for every local
sentence $\varphi$, whenever $\varphi$ has a model of order type
$\mu$, then $\varphi$ has models of order type $\tau$, for each
cardinal $\tau\geq \om$.
\end{Deff}

\noi It is proved in \cite{fr}, using the sentence $\phi_0$ given by Theorem \ref{phi0},
 that $\mu$ is at least inaccessible. We can now state the folllowing result.

\begin{The}
 The ordinal $\mu$ is  the first $\om$-Mahlo cardinal, if such a cardinal exists.
\end{The}

\proo
Ressayre proved in \cite{ress} that $\mu$ exists if an $\om$-Mahlo cardinal exists and then $\mu$ is bounded
by this large cardinal.
On the other hand, using the sentence $\phi_n$ given by Theorem \ref{phin},
we can see that $\mu$ must be an $n$-Mahlo cardinal.  The ordinal $\mu$ is then an $\om$-Mahlo cardinal because it is
an $n$-Mahlo cardinal, for each integer $n\geq 1$.
Thus the ordinal $\mu$ is in fact the first $\om$-Mahlo cardinal, if such a cardinal exists.

\section{Strength of the stretching principles}

We are going to prove now the following result.

\begin{The}
The statement $(E)$  (respectively, the statement $(E')$) implies
the existence, for each integer $n$,
 of an $n$-Mahlo cardinal.
\end{The}

\proo
 Assume first that $(E)$ is true in a model $U$ of ZF.  We know from  \cite{fr} that there exists a local
sentence $\varphi$ such that $Sp_\infty(\varphi)=\{\om\}$. Consider
now the sentence $\Psi_n=\varphi \cup \phi_n$, where $\phi_n$ is
given by Theorem \ref{phin}. Assume that there is no $n$-Mahlo
cardinal in $U$. Then the  sentence  $\Psi_n$ has models of order
type $\kappa$, for each cardinal $\kappa\geq \om$.  So if  $(E)$ is
true, then $\Psi_n$ would have a finite model $M$  containing
$N_{\Psi_n}$ remarkable and monotonic indiscernibles. But $M$ can
not satisfy $\varphi$ because otherwise the stretching $M(\om +1)$
would be a well ordered model of    $\varphi$  of order type greater
than $\om$. (recall that by Lemma 8 of \cite{fr},if the
indiscernibles are monotonic, then, for each ordinal $\alpha$, the
stretching $M(\alpha)$ is  well ordered) . And $M$ can not satisfy
$\phi_n$ because otherwise the sentence $\phi_n$ would have a model
of order type $\om$ by the equivalence $(C)$ of the Stretching
Theorem \ref{stretching}. This would lead to a contradiction so we
can conclude that $U$ contains some $n$-Mahlo cardinal for each
integer $n\geq 1$. \nl To prove the corresponding result for $(E')$,
it suffices to add a unary predicate $R$ to the signatures of
$\varphi$ and $\phi_n$ and to reason in a similar way, replacing the
sentence $\varphi$ by the sentence $\varphi \wedge \fa x R(x)$ and
the sentence $\phi_n$ by the sentence $\phi_n \wedge \fa x R(x)$.

\hs We recall from \cite{fr} that  $(E)$  and  $(E')$ can be proved assuming the existence for  each integer $n$
 of an $n$-Mahlo cardinal. Thus we can infer the following result.

\begin{Cor}
The statement $(E)$  (respectively, the statement $(E')$) is
equivalent to the existence, for each integer $n$,
 of an $n$-Mahlo cardinal.
\end{Cor}

\hs

\end{document}